\documentclass[a4paper]{article}
\usepackage{amsfonts}
\usepackage{amsmath}
\usepackage{amssymb}
\begin{document}
\newtheorem{theorem}{Theorem}[section]
\newtheorem{lemma}[theorem]{Lemma}
\newtheorem{corollary}[theorem]{Corollary}
\newtheorem{conjecture}[theorem]{Conjecture}
\newtheorem{remark}[theorem]{Remark}
\newtheorem{definition}[theorem]{Definition}
\newtheorem{problem}[theorem]{Problem}
\newtheorem{example}[theorem]{Example}
\newtheorem{proposition}[theorem]{Proposition}
\title{{\bf Curvature semipositivity of relative pluricanonical systems}}
\date{March 28, 2007}
\author{Hajime TSUJI\footnote{This work was supported by JSPS fellowship.}}
\maketitle
\begin{abstract}\noindent We prove semipositivity of the curvature of the Narashimhan-Simha metric on 
an arbitrary flat projective family  of varieties with only canonical singularities. This also implies that the direct image of pluricanonical systems 
have natural continuous hermitian strucutres with semipositive curvature curvature currents in the sense of Nakano.  This paper is a revised version of  \cite{tu6}
using \cite{tu7,tu8}. 
\end{abstract}
\tableofcontents
\section{Introduction}

Let $f : X \longrightarrow S$ be a flat projective family of varieties over a complex manifold $S$.
In \cite{n-s}, Narashimhan and Simha constructed  a continuous  hermitian metric $h_{m}$ on  $mK_{X/S}$  for every sufficiently large positive integer $m$, when $f$ is smooth and  $K_{X/S}$ is relatively ample. 
They studied the parameter dependence of the metrics 
(\cite[p.122, Section 3]{n-s}) 
and proved the existence of moduli space of smooth projective varieties 
with ample canonical bundle (\cite[p.126, Theorem]{n-s}). 

This metric depends only on the complex structure of $X$, hence it is functorial under biholomorphisms. 
Although they argued only the case of relatively ample $K_{X/S}$, 
the same construction can be applied 
also to the case that a general fiber has nonnegative Kodaira dimension (cf. Section 3). 
We call the metric the Narashimhan-Simha metric. 
But in this case we see that the resulting metric is a singular hermitian metric.  And the metric is functorial under birational morphism in an appropriate sense (cf.
Section \ref{invariance}).   
  
 The semipositivity result of Y. Kawamata \cite[p.57, Theorem 1]{ka1} 
 proved that when $X$ is smooth and $S$ is a projective curve, then $f_{*}{\cal O}_{X}(mK_{X/S})$ is semipositive on $S$ in the sense of algebraic geometry, i.e., any quotient sheaf of $f_{*}{\cal O}_{X}(mK_{X/S})$ has semipositive degree.
We note that the semipositivity  in \cite[p.66,Theorem 1]{ka1} comes from a 
Finsler type metric but not from a hermitian metric on $f_{*}{\cal O}_{X}(mK_{X/S})$.

The purpose of this paper is to refine this semipositivity result 
in terms of the Narashimhan-Simha singular hemitian metric on the multi 
relative canonical sheaf ${\cal O}_{X}(mK_{X/S})$.  
In this case it is natural to consider flat projective families of varieties with only canonical singularities (cf. Definition \ref{Can Sing}).  
Namely we prove the following theorem. 

\begin{theorem}\label{main theorem}(\cite[Theorem 1.1]{tu6})
Let $f : X \longrightarrow S$ be a flat projective family with connected fibers  
such that a general fiber has only canonical singularities over a complex manifold $S$.  We set  
\[
S^{\circ} := \{ s\in S\mid \mbox{$X_{s} := f^{-1}(s)$ is reduced and has only canonical singularities}\}.
\]
Let $h_{m}$ be the Narashimhan-Simha singular hermitian metric  on $mK_{X/S}$
(cf. Section \ref{N-S}) on $X^{\circ} := f^{-1}(S^{\circ})$
 
Then $h_{m}$ has semipositive curvature in the sense of current on $X^{\circ}$ 
and the curvature $-dd^{c}\log h_{m}$ extends across $X - X^{\circ}$ 
as a closed positive current. In particular $h_{m}$ extends to a singular hermitian metric on $mK_{X/S}$ across $X - X^{\circ}$. $\square$
\end{theorem}
The main difference between Theorem 1.1 and  \cite[p.57, Theorem 1]{ka1}
is that in Theorem 1.1 the positivity is on the family $X$ not on  the base space $S$.
This fact plays the key role in this paper. Also we note that it is clear that 
 $h_{m}$ has semipositive curvature in the vertical direction, i.e., the restriction of $h_{m}$ to every fiber of $f$ has semipositive curvature.  Hence the main assertion of Theorem \ref{main theorem} is that  the curvature of $h_{m}$ is semipositive also in the horizontal direction.  
This theorem is new even when $K_{X/S}$ is relatively ample 
and $f : X \longrightarrow S$ is smooth.  \vspace{2mm} \\
We note that the statement of Theorem 1.1 is meaningful, only if 
$f_{*}{\cal O}_{X}(mK_{X/S})$ is not zero. 
Otherwise $h_{m}$ is not defined (cf. Section \ref{N-S}).

As a direct consequence of Theorem \ref{main theorem}, we have 
the following theorem.
 
\begin{theorem}\label{semipositive}
Let $f : X \longrightarrow S$ be a flat projective family with connected fibers  such that a general fiber of $f$ has only canonical singularities. 
We assume $X$ is normal and $S$ is smooth. 
We set  
\[
S^{\circ} := \{ s\in S\mid \mbox{$X_{s} := f^{-1}(s)$ is reduced and has only canonical singularities}\}.
\] 
Then $K_{X/S}$ has a relative AZD $h$ over $S^{\circ}$ such that 
$\Theta_{h}$ is semipositive on $X$ (Theorem \ref{family}). 

And $F_{m} := f_{*}{\cal O}_{X}(mK_{X/S})$ carries a  hermitian metric $h_{F_{m}}$ defined by 
\[
h_{F_{m}}(\sigma ,\sigma^{\prime})_{s}:= 
(\sqrt{-1})^{n^{2}}\int_{X_{s}}\sigma\wedge \bar{\sigma}^{\prime}\cdot h^{m-1}
\hspace{10mm}(\sigma ,\sigma^{\prime}\in F_{m,s},s\in S^{\circ})
\]
with Nakano semipositive curvature in the sense of current over
$S^{\circ}$, where $n = \dim X - \dim S$. 
$\square$
\end{theorem}
\begin{remark}\label{fg}
According to the recent outstanding result on finite generation of canonical rings (\cite{b-c-h-m}), one may use $h_{m}^{1/m}$ instead of $h$ for the relative AZD for some positive 
integer $m$.  By the definition $h_{m}^{1/m}$ (see Section 3.1) have only algebraic singularities (cf. Definition \ref{alg sing}).  In this case the resulting metric on $F_{m}$ seems to be quite regular. This might be a great 
advantage for the future applications.  $\square$
\end{remark}

\noindent We note that there is another approach to obtain Theorem \ref{semipositive}
using the dynamical construction of AZD's (cf \cite{tu8}).   

This article  consists of three parts.
The first part (Section 1,2) consists of  the introduction and preliminaries. 
The second part (Section 3,4) consists of the proof of the semipositivity theorem, Theorem 1.1 and the construction of the canonical AZD for the relative canonical bundle.  The third part (Section 5) consists of the proof of 
Theorem \ref{semipositive}. 

The author would like to express scincere gratitude to Professor Mihai Paun 
who pointed out an error in the proof of Theorem \ref{nakano} and the reference
\cite{b-p}. 
  
In this paper all the varieties are defined over $\mathbb{C}$ and the Bergman
kernels means the diagonal parts (we do not use the reproducing kernel properties in this paper). 

\section{Preliminaries}
Here we shall collect basic notions which are used in this article. 

\subsection{Canonical singularities}

The notion of canonical singularities were introduced by M. Reid to study 
canonical models of general type (\cite{r}). 
Actually canonical singularities are caractertized as singularities which 
appear on the canonical models when the canonical rings are finitely generated.

\begin{definition}\label{Can Sing}
Let $X$ be a normal variety such that the canonical divisor $K_{X}$ is 
{\bf Q}-Caritier. 
If $\mu : Y \longrightarrow X$ is a resolution of $X$, then we can write
\[
K_{Y} = \mu^{*}K_{X} + F
\]
with $F = \sum_{j}e_{j}E_{j}$ for the exceptional divisors $\{ E_{j}\}$. 
We call $F$ the discrepancy and $e_{j}\in  \mbox{\bf Q}$ the discrepancy
coefficient for $E_{j}$. 

$X$ is said to have only {\bf terminal singularities}
(resp. {\bf canonical singularities}, if 
$e_{j} > -1$ (resp. $\geqq -1$)
for all $j$ for a log resolution $\mu : Y \longrightarrow X$.
One can also say that $X$ is terminal (resp. canonical), or $K_{X}$ is terminal (resp. canonical), when $X$ has only terminal (resp. canonical) singularities.
$\square$ \end{definition}

\begin{definition}\label{index}
Let $X$ be a normal variety with only canonical singularities. 
Let $x$ be a point on $X$. 
The {\bf local index} of $X$ at $x$ is the least positive integer $r$ such that 
$rK_{X}$ is Cartier near $x$. 
The {\bf global index} of $X$ is the least positive integer $r$ such that 
$rK_{X}$ is Cartier on $X$. $\square$
\end{definition}
The following theorem is fundamental in this article.

\begin{theorem}(\cite[p.85,Main Theorem]{ka2})\label{defcan}
Let $\pi : X \longrightarrow S$ be a flat morphism from a germ of variety 
$(X,x_{0})$ to a germ of a smooth curve $(S,s_{0})$ whose spetial fiber 
$X_{0} = \pi^{-1}(s_{0})$ has only canonical singularities.
Then $X$ has only canonical singularities.
In particular, the fibers $X_{s} = \pi^{-1}(s)$ have only canonical singularities. $\square$
\end{theorem}

\subsection{Singular hermitian metrics}\label{singh}
In this subsection $L$ will denote a holomorphic line bundle on a complex manifold $M$. 
\begin{definition}\label{singhm}
A  {\bf singular hermitian metric} $h$ on $L$ is given by
\[
h = e^{-\varphi}\cdot h_{0},
\]
where $h_{0}$ is a $C^{\infty}$-hermitian metric on $L$ and 
$\varphi\in L^{1}_{loc}(M)$ is an arbitrary function on $M$.
We call $\varphi$ a  weight function of $h$. $\square$ 
\end{definition}
The curvature current $\Theta_{h}$ of the singular hermitian line
bundle $(L,h)$ is defined by
\[
\Theta_{h} := \Theta_{h_{0}} + \sqrt{-1}\partial\bar{\partial}\varphi ,
\]
where $\partial\bar{\partial}$ is taken in the sense of a current.
The $L^{2}$-sheaf ${\cal L}^{2}(L,h)$ of the singular hermitian
line bundle $(L,h)$ is defined by
\[
{\cal L}^{2}(L,h) := \{ \sigma\in\Gamma (U,{\cal O}_{M}(L))\mid 
\, h(\sigma ,\sigma )\in L^{1}_{loc}(U)\} ,
\]
where $U$ runs over the  open subsets of $M$.
In this case there exists an ideal sheaf ${\cal I}(h)$ such that
\[
{\cal L}^{2}(L,h) = {\cal O}_{M}(L)\otimes {\cal I}(h)
\]
holds.  We call ${\cal I}(h)$ the {\bf multiplier ideal sheaf} of $(L,h)$.
If we write $h$ as 
\[
h = e^{-\varphi}\cdot h_{0},
\]
where $h_{0}$ is a $C^{\infty}$ hermitian metric on $L$ and 
$\varphi\in L^{1}_{loc}(M)$ is the weight function, we see that
\[
{\cal I}(h) = {\cal L}^{2}({\cal O}_{M},e^{-\varphi})
\]
holds.
For $\varphi\in L^{1}_{loc}(M)$ we define the multiplier ideal sheaf of $\varphi$ by 
\[
{\cal I}(\varphi ) := {\cal L}^{2}({\cal O}_{M},e^{-\varphi}).
\] 
\begin{example}
Let $\sigma\in \Gamma (X,{\cal O}_{X}(L))$ be the global section. 
Then 
\[
h := \frac{1}{\mid\sigma\mid^{2}} = \frac{h_{0}}{h_{0}(\sigma ,\sigma)}
\]
is a singular hemitian metric on $L$, 
where $h_{0}$ is an arbitrary $C^{\infty}$-hermitian metric on $L$
(the right hand side is ovbiously independent of $h_{0}$).
The curvature $\Theta_{h}$ is given by
\[
\Theta_{h} = 2\pi\sqrt{-1}(\sigma )
\]
where $(\sigma )$ denotes the current of integration over the 
divisor of $\sigma$. $\square$ 
\end{example}
\begin{definition}\label{alg sing}
If $\{\sigma_{i}\}$ are a finite number of global holomorphic sections of $L$,
for every positive rational number $\alpha$ and a continuous function 
$\phi$,
\[
h := e^{-\phi}\cdot\frac{1}{\sum_{i}\mid\sigma_{i}\mid^{2\alpha}}
\]
defines a singular hermitian metric  on 
$\alpha L$.
We call such a metric $h$ a singular hermitian metric 
on $\alpha L$ with  {\bf algebraic singularities}. $\square$
\end{definition}
Singular hermitian metrics with algebraic singularities 
are particulary easy to handle, because its multiplier 
ideal sheaf of the metric can 
be controlled by taking  suitable successive blowing ups 
such that the total transform of the divisor
$\sum_{i}(\sigma_{i})$ is a divisor with normal crossings. 
\begin{definition}\label{pe}
$L$ is said to be {\bf pseudoeffective}, if there exists 
a singular hermitian metric $h$ on $L$ such that 
the curvature current 
$\Theta_{h}$ is a closed positive current.
Also a singular hermitian line bundle $(L,h)$ is said to be {\bf pseudoeffective}, 
if the curvature current $\Theta_{h}$ is a closed positive current. $\square$
\end{definition}

\subsection{Lelong numbers}\label{Lelong number}

The Lelong number is a measure of singularities of closed positive 
currents.   
We note that a subvariety of codimension $p$ can be identified with 
the current of integration over it.  
In this special case the Lelong number exactly corresponds to the multiplicity.
Hence we can consider the Lelong number as a generalization of the notion
of  multiplicities.  \vspace{3mm} \\   

Let $W\subset \mathbb C^n$ be a domain, and $\Theta$ a positive
current of degree $(q,q)$ on $W$. For a point $p\in W$ one defines
$$
\nu(\Theta,p,r)=\frac{1}{r^{2(n-q)}}\int_{\|z-p\|<r} \Theta(z)
\wedge (dd^c\|z\|^2)^{n-q}.
$$
The {\bf Lelong number} of $\Theta$ at $p$ is defined as
$$
\nu(\Theta,p)= \lim_{{r\to 0}\atop{r\; >\; 0}}\nu(\Theta,p,r).
$$
If $\Theta$ is the curvature of $h= e^{-u}$, $u$
plurisubharmonic, one has
\begin{equation}\label{liminf}
\nu(\Theta,p)= \sup \{\gamma \geq 0; u\leqq \gamma\log(\|z-p\|^2) +
O(1)\}.
\end{equation}
The definition of a singular hermitian metric carries over to the
situation of reduced complex spaces.

\begin{definition}\label{singdef} Let $Z$ be a reduced complex space
and $L$ a holomorphic line bundle. A {\em singular} hermitian
metric $h$ on $L$ is a singular hermitian metric $h$ on
$L|Z_{reg}$ with the following property: There exists a
desingularization $\pi: \tilde{Z} \to Z$ such that $h$ can be
extended from $Z_{reg}$ to a singular hermitian metric $\tilde{h}$
on $\pi^*L$ over $\tilde{Z}$. $\square$
\end{definition}

The definition is independent of the choice of a desingularization
under a further assumption. Suppose that $\Theta_{\tilde{h}} \geq -c
\cdot \omega $ in the sense of currents, where $c>0$, and $\omega$
is a positive definite, real $(1,1)$-form on $\tilde{Z}$ of class
$C^{\infty}$. Let $ \pi_1 :{Z_1} \to Z$ be a further
desingularization. Then $\tilde{Z}\times_Z Z_1 \to Z$ is dominated
by a desingularization $Z'$ with projections $p :Z' \to \tilde{Z} $
and $p_1 :Z' \to Z_1$. Now $p^*\log \tilde{h}$ is of class
$L^1_{loc}$ on $Z'$ with a similar lower estimate for the
curvature. The push-forward $p_{1*}p^*\tilde{h}$ is a singular
hermitian metric on $Z_1$. In particular, the extension of $h$ to
a desingularization of $Z$ is unique.

In \cite{g-r} for plurisubharmonic functions on a normal complex
space the Riemann extension theorems were proved, which will be
essential for our application. 

For a reduced complex space a  plurisubharmonic function $u$ is by
definition an upper semi-continuous function $u:X \to [-\infty,
\infty)$ whose restriction to any local, smoothly parameterized
analytic curve is either identically $-\infty$ or subharmonic.

A locally bounded function $u:X \to [-\infty, \infty)$ from
$L^1_{loc}(X)$ is called weakly plurisubharmonic, if its restriction to
the regular part of $X$ is plurisubharmonic.

Differential forms with compact support on a reduced complex
space are by definition locally extendable to an ambient
subspace, which is an open subset $U$ of some $\mathbb C^n$. Hence
the dual spaces of differential $C^\infty$-forms on such $U$
define currents on analytic subsets of $U$. The positivity of a
real $(1,1)$-current is defined in a similar way as above.

For locally bounded functions of class $L^1_{loc}$ the weak
plurisubharmonicity is equivalent to the positivity of the
current $d d^c u$. It was shown that these functions are exactly
those, whose pull-back to the normalization of $X$ are plurisubharmonic. We
note

\begin{definition}
Let $L$ be a holomorphic line bundle on a reduced complex space
$X$. Then a singular hermitian metric $h$ is called positive, if
the functions, which define locally $-\log h$ are weakly plurisubharmonic.
$\square$
\end{definition}

This definition is compatible with Definition \ref{singhm}. Let
$L$ be a holomorphic line bundle on a complex space $Z$ equipped
with a positive, singular hermitian metric $h_r$ on $L|Z_{reg}$.
If $\pi:\tilde Z \to Z$ is a desingularization, and $\tilde h$ a
positive, singular hermitian metric on $\pi^* L$, extending
$h|Z_{reg}$, we see that $-\log h_r$ is locally bounded at the
singularities of $Z$ so that $\tilde{h}$ induces a singular, positive
metric on $L$ over $Z$. \vspace{5mm} \\
The  Lelong number of the curvature current of a pseudoeffective
singular hemitian line bundle (cf. Definition \ref{pe}) can be observed as an obstruction for 
nefness by the following lemma. 

\begin{lemma}(cf. \cite{tu,tu2})\label{nefness}
Let $X$ be a smooth projective variety and let $(L,h)$ be a singular hermitian 
line bundle such that $\Theta_{h}$ is semipositive in the sense of current 
and the Lelong number $\nu(\Theta_{h})$ is identically $0$ on $X$.
Then $L$ is nef. 
\end{lemma}
{\bf Proof of Lemma \ref{nefness}}.
Let $x\in X$ be an arbitrary point on $X$.  
Let $A$ be a positive line bundle on $X$ and let $h_{A}$ be a $C^{\infty}$-hermitian metric on $A$ such that $\omega := \Theta_{h_{A}}$ is a K\"{a}hler form.
And let $d_{x}$ denote the distance
function from $x$ with respect to a K\"{a}hler form $\omega$ on $X$.
Let $a$ be a positive integer such that 
\[
(2n+2)\log d_{x} + \mbox{Ric}_{\omega} +a\Theta_{h_{A}} 
\]
is strictly positive for every $x\in X$, where $n$ denotes the dimension of $X$.
This implies that
\[
d_{x}^{-(2n+2)}\cdot h_{A}^{a} \cdot (\omega^{n})\cdot h_{L}^{m}
\] 
is a singular hermitian metric 
on $A - K_{X} + mL$ with strictly positive curvature. 
Then by Nadel's vanishing theorem \cite[p.561]{n} and the description of 
the Lelong number (\ref{liminf}), we see that $mL + aA$
is free at $x$ for every $m$. 
Since $m$ and $x$ are arbitrary, this implies that $L$ is nef. $\square$
\begin{remark}
It is easy to generalize Lemma \ref{nefness}, in the form : \\
If $(L,h)$ is a pseudoeffective singular hermitian line bundle on a 
smooth projective variety $X$.  Then for every $x \in X$ such that $\nu (\Theta_{h},x) = 0$ and an irreducible curve $C$ containing $x$, 
$L\cdot C \geqq 0$ holds, i.e., $L$ is nef modulo the positive Lelong number 
locus of $\Theta_{h}$. $\square$
\end{remark}

\subsection{Analytic Zariski decompositions (AZD)}
In this subsection we shall introduce the notion of analytic Zariski decompositions. 
By using analytic Zariski decompositions, we can handle  big line bundles
like  nef and big line bundles.
\begin{definition}\label{defAZD}
Let $M$ be a compact complex manifold and let $L$ be a holomorphic line bundle
on $M$.  A singular hermitian metric $h$ on $L$ is said to be 
an analytic Zariski decomposition, if the followings hold.
\begin{enumerate}
\item $\Theta_{h}$ is a closed positive current,
\item for every $m\geq 0$, the natural inclusion
\[
H^{0}(M,{\cal O}_{M}(mL)\otimes{\cal I}(h^{m}))\rightarrow
H^{0}(M,{\cal O}_{M}(mL))
\]
is an isomorphim. $\square$
\end{enumerate}
\end{definition}
\begin{remark} If an AZD exists on a line bundle $L$ on a smooth projective
variety $M$, $L$ is pseudoeffective by the condition 1 above. $\square$
\end{remark}

\begin{theorem}(\cite{tu,tu2})
 Let $L$ be a big line  bundle on a smooth projective variety
$M$.  Then $L$ has an AZD. $\square$
\end{theorem}
As for the existence for general pseudoeffective line bundles, 
now we have the following theorem.
\begin{theorem}(\cite[Theorem 1.5]{d-p-s})\label{AZD}
Let $X$ be a smooth projective variety and let $L$ be a pseudoeffective 
line bundle on $X$.  Then $L$ has an AZD.$\square$
\end{theorem}

\section{Narashimhan-Simha metrics} 

In this section, following \cite{n-s}, we shall define a singular hermitian metric on 
the canonical bundle of a smooth projective manifold with nonnegative Kodaira dimension and study the behavior of the metric under deformation.

\subsection{Narashimhan-Simha metrics}\label{N-S}
Let $X$ be a smooth projective variety with nonnegative Kodaira dimension of dimension $n$. 
Let $m$ be a positive integer.
For a section $\eta \in H^{0}(X,{\cal O}_{X}(mK_{X}))$ we define  
a nonnegative number $\parallel\eta\parallel_{\frac{1}{m}}$ by 
\[
\parallel\eta\parallel_{\frac{1}{m}} = \mid \int_{X}(\eta\wedge\bar{\eta})^{\frac{1}{m}}\mid^{\frac{m}{2}}.
\]
Then $\eta\mapsto \parallel\eta\parallel_{\frac{1}{m}}$ is a continuous pseudonorm 
on $H^{0}(X,{\cal O}_{X}(mK_{X}))$, i.e.,  it is a continuous and has the properties :
\begin{enumerate}
\item $\parallel\eta\parallel_{\frac{1}{m}} = 0 \Leftrightarrow \eta = 0$,
\item $\parallel\lambda \eta\parallel_{\frac{1}{m}} = \mid\lambda\mid\cdot
\parallel\eta\parallel_{\frac{1}{m}} $ holds for all $\lambda\in \mathbb{C}$. 
\end{enumerate}
But it is not a norm on $H^{0}(X,{\cal O}_{X}(mK_{X}))$ except $m=1$.
We define a continuous section $K_{m}$ of 
$(K_{X}\otimes \bar{K}_{X})^{\otimes m}$ 
\[
K_{m}(x):= \sup\{ \mid \eta (x)\mid^{2}\,\,\,; 
\,\,\,\parallel \eta \parallel_{\frac{1}{m}} = 1\}  \,\,\,\,\,\,\, (x\in X),  
\]
where $\mid\eta (x)\mid^{2} = \eta (x)\otimes \bar{\eta}(x)$.
We call $K_{m}$ {\bf the $m$-th Narashimhan-Simha potential} of $X$.
The continuity of $K_{m}$ follows from the expression 
\[
K_{m}(x) = \sup \{ \mid \eta (x)\mid^{2}/\parallel\eta\parallel_{\frac{1}{m}}\,\, ; \,\,\,\,
\eta = \sum_{i=0}^{N(m)} c_{i}\sigma_{i}^{(m)}, \sum_{i=0}^{N(m)}\mid c_{i}\mid^{2} = 1\} ,
\]
where $\{ \sigma_{0}^{(m)},\cdots ,\sigma_{N(m)}^{(m)}\}$ 
is a basis for $H^{0}(X,{\cal O}_{X}(mK_{X}))$. 
We define the singular hermitian metric $h_{m}$ on 
$H^{0}(X,{\cal O}_{X}(mK_{X}))$ by 
\[
h_{m} := \frac{1}{K_{m}}.
\]
We call $h_{m}$ {\bf the Narashimhan-Simha metric}										 on $mK_{X}$. 
This metric is introduced by Narashimhan and Seshadri for smooth
 canonically polarized variety to study the moduli space of canonically polarized variety (\cite{n-s}).  
We note that the singularities of $h_{m}$ is located exactly on the support of 
the base locus of $\mid mK_{X}\mid$. 
Hence in the case of canonically polarized variety, $h_{m}$ is a nonsingular continuous metric on $mK_{X}$ for every sufficiently large $m$.  
Since $K_{m}$ is locally the supremum of the  family of powers of absolute value of  holomorphic functions we see that the curvature
\[
\Theta_{h_{m}} := \sqrt{-1}\partial\bar{\partial}\log K_{m}
\] 
is a closed positive current. 

For applications sometimes more convenient to consider the following 
regularization of $\Theta_{h_{m}}$.
  
Let us define the $L^{2}$-inner product on $H^{0}(X,{\cal O}_{X}(mK_{X}))$ by 
\[
(\eta ,\tau)_{m} := (\sqrt{-1})^{n^{2}}\int_{X}\frac{\eta\cdot\bar{\tau}}{K_{m}^{\frac{m-1}{m}}}
\hspace{10mm} (\eta ,\tau \in H^{0}(X,{\cal O}_{X}(mK_{X}))). 
\]
The inner product is well defined, because by the definition of 
$K_{m}$, the zero of $K_{m}$ corresponds exactly the ideal of 
$\mbox{Bs}\mid\!\!mK_{X}\!\!\mid$.
Let $\{ \sigma_{0}^{(m)},\cdots \sigma_{N(m)}^{(m)}\}$ be a complete orthonormalbasis of $H^{0}(X,{\cal O}_{X}(mK_{X}))$ with respect to the innner product
$(\,\, ,\,\, )_{m}$. 
We define the Bergman kernel $\tilde{K}_{m}$ by 
\[
\tilde{K}_{m} = \sum_{i=0}^{N(m)}\mid\!\sigma_{i}^{(m)}\!\!\mid^{2}.
\]
Then $\tilde{K}_{m}$ is independent of the choice of the complete orthonormal basis $\{ \sigma_{0}^{(m)},\cdots \sigma_{N(m)}^{(m)}\}$.
In fact 
\[
\tilde{K}_{m}(x) := \sup \{ \mid\sigma (x)\mid^{2}\mid 
\sigma \in H^{0}(X,{\cal O}_{X}(mK_{X})), \parallel\sigma\parallel = 1\} 
\]
holds, 
where $\parallel\,\,\,\,\parallel$ denotes the $L^{2}$-norm with respect to 
$(\,\,\,\, ,\,\,\,\,)_{m}$ (see for example \cite[p.46]{kr}).

Then the closed positive current
\[
\omega_{m} := \sqrt{-1}\partial\bar{\partial}\log \tilde{K}_{m}
\]
coincides (up to some normalized constant) the pullback of the Fubini-Study K\"{a}hler form on $\mathbb{P}^{N(m)}$ 
on $X \backslash \mbox{Supp}\,\mbox{Bs}\mid\!mK_{X}\!\!\mid$ 
by the rational map : 
\[
X\ni x -\cdots\rightarrow [\sigma_{0}^{(m)}(x):\cdots :\sigma_{N(m)}^{(m)}(x)]
\in \mathbb{P}^{N(m)}.
\] 
We call $\omega_{m}$ {\bf the Narashimhan-Shimha K\"{a}hler current}.

\subsection{Dependence of Narashimhan-Simha metrics on families}\label{continuity}

Let $f : X \longrightarrow S$ be a flat projective family of varieties with only canonical singularities. 
Let $n$ denote the relative dimension of $f : X \longrightarrow S$. 
Let $X_{s}$ denote the fiber $f^{-1}(s) (s\in S)$.  
Let $m$ be a positive integer and let $K_{m,s}$ be the Narashimhan-Shimha 
potential of $mK_{X_{s}}$.  
Then we obtain the family $\{ mK_{X_{s}}\}_{s\in S}$ of Narashimhan-Simha 
potentials over $S$ and the family of Narashimhan-Simha metrics 
$\{ h_{m,s}\}_{s\in S}$ on $X_{s}$.

In \cite{n-s} Narashimhan-Simha proved the continuous dependence of 
$K_{m,s}$ with respect to $s\in S$, when $K_{X/S}$ is relatively ample and $f$ is smooth. 
The key argument here is the invariance of the plurigenera $P_{m}(X_{s})
= h^{0}(X_{s},{\cal O}_{X_{s}}(mK_{X_{s}}))$ with respect to $s\in S$ which is derived by the Kodaira vanishing theorem at that time.  

We note that by \cite[p.91,Theorem 6]{ka2} and \cite{s1,s2,tu3}, the plurigenera $P_{m}(X_{s})$ is independent of 
$S$.  Hence every global holomorphic section of $mK_{X_{s}}$ extends holomorphically to nearby fibers.  
  
Hence by the defiition of the Narashimhan-Simha potential, we have the following lemma.  

\begin{lemma}\label{conti}
$\{ K_{m,s}\}$ is continuous with on $X$. $\square$
\end{lemma}

By Lemma \ref{conti} we see that there exists a relative volume form 
$K_{m}$ such that 
\[
K_{m}\mid_{X_{s}} = K_{m,s}
\]
holds for all $s\in S$.
We call $K_{m}$ the {Narashimhan-Simha potential} of 
the family $f : X \longrightarrow S$ and also 
\[
h_{m} := \frac{1}{K_{m}}
\]
the {\bf Narashimhan-Simha} metric on $mK_{X/S}$ respectively. 

\subsection{Positivity of Narashimhan-Simha metrics on projective families}\label{positive}

In this subsection we shall prove Theorem 1.1. 
Let $f : X \longrightarrow S$ be a flat projective family with connected fibers  
such that a general fiber has only canonical singularities over a complex manifold $S$.  We set  
\[
S^{\circ} := \{ s\in S\mid \mbox{$X_{s} := f^{-1}(s)$ is reduced and has only canonical singularities}\}.
\]
For $s\in S$, we denote $X_{s}$ the fiber $f^{-1}(s)$. 
For the proof we may and do assume that $S$ is the open unit disk
$\Delta$ in $\mathbb{C}$.   We note that $f^{*}{\cal O}_{X}(mK_{X/S})$
is not compatible with base change.   
Nevertheless this is enough to prove the statement of Theorem \ref{main theorem} on 
the extension of $h_{m}$ across $X - X-^{\circ}$ holds, by slicing. 

Let us denote $f_{*}{\cal O}_{X}(mK_{X/S})$ by $F_{m}$.

\begin{lemma}(\cite[p.63, Lemma 7 and p.64, Lemma 8]{ka1})\label{subh}
Let $\eta \in H^{0}(S,{\cal O}_{S}(F_{m}))$ be a holomorphic section. 
Let  $\parallel\!\eta (s)\!\parallel_{\frac{1}{m}}$ be the Narashimhan-Sihma pseudonorm of $\eta (s) \in H^{0}(X_{s},{\cal O}_{X_{s}}(mK_{X_{s}}))$.
Then 
\[
\log \parallel\!\eta (s)\!\parallel_{\frac{1}{m}} \geqq \frac{1}{2\pi}\int_{0}^{2\pi}\log \parallel\!\eta (s +re^{\sqrt{-1}\theta})\!\parallel_{\frac{1}{m}}d\theta 
\]
holds for every $s\in S$ and $r > 0$ such that $\Delta (s,r)\subset S$, i.e.
$\log \parallel\!\eta (s)\!\parallel_{\frac{1}{m}}$ is superharmonic on $S$, where $\Delta (s,r)$ denotes the open disk of radius $r$ with centre $s$.  $\square$
\end{lemma}
\begin{remark}
In \cite{ka1}, Lemma \ref{subh} was considered only for smooth projective families.  But by Theorem \ref{defcan} the same proof works for flat projective families with only canonical singularities.  We note that canonical singularities are rational (\cite{e} or cf.\cite{k-m}). $\square$
\end{remark}
Let us fix an arbitrary $s\in S^{\circ}$.  
Let $x\in X_{s}$, $r$ be a positive number such that 
$\Delta (s,r) \subset S$ and let 
$\tau : \Delta (s,r) \longrightarrow X$ be any holomorphic section 
such that $\tau (s) = x$. 
In general the such a section does not exist when $x$ is a singular point on 
$X_{s}$.  This case will be treated later. 

Let $\eta \in H^{0}(S,{\cal O}_{S}(F_{m}))$ such that 
$\parallel \eta (s)\parallel_{\frac{1}{m},s} = 1$ and 
\[
\mid\eta(s,x)\mid^{2} = K_{m,s}(x)
\]
holds, where $\eta (s)$ denotes the restriction $\eta\mid_{X_{s}}$ and $\eta (s,x)$ denotes the point value of $\eta (s)$ at $x$.

Let $a$ be a positive integer such that $aK_{X/S}$ is Cartier. 
Let $\Omega$ be a local  generator of the invertible sheaf ${\cal O}_{X}(aK_{X/S})$ on a neighbourhood $U$ of 
$x$.   Shrinking  $S$, if necessary, we may assume that $\tau (\Delta (s,r))\subset U$
holds. 
Let $f$ be the function on $\Delta (s,r)$ defined by 
\[
f(z) = \mid\frac{\eta}{\Omega^{m/a}}\mid^{2}(\tau (z)) \,\,\,\,\,\, (z\in \Delta (s,r)) 
\]
Then we see that 
\begin{equation}
\log f(s) \leqq \frac{1}{2\pi}\int_{0}^{2\pi}\log f(s+re^{\sqrt{-1}\theta}) d\theta
\end{equation}
holds by the subharmonicity of the logarithm of abolute value of holomorphic function. 
On the other hand we see that 
\begin{equation}
\log \parallel\eta (s)\parallel_{\frac{1}{m}} \geqq 
\frac{1}{2\pi}\int_{0}^{2\pi}\log \parallel\eta(s+re^{\sqrt{-1}\theta})
\parallel_{\frac{1}{m}}d\theta 
\end{equation}
holds by Lemma \ref{subh}. 

By (1) and (2), we see that 
\[
\log \frac{f(s)}{\parallel\eta (s)\parallel_{\frac{1}{m}}}\leqq \frac{1}{2\pi}\int_{0}^{2\pi}\log \frac{f(s+re^{\sqrt{-1}\theta})}{\parallel\eta (s+re^{\sqrt{-1}\theta}))\parallel_{\frac{1}{m}}}d\theta
\]
We note that 
\[
K_{m}(\tau (s+re^{\sqrt{-1}\theta})) \geqq 
\frac{f(s+re^{\sqrt{-1}\theta})}{\parallel\eta (s+re^{\sqrt{-1}\theta}))\parallel_{\frac{1}{m}}}\mid \Omega\mid^{\frac{m}{a}}
\]
holds by the definition of $K_{m}$. 
Hence we see that 
\[
\log K_{m}(\tau (s)) \leqq \frac{1}{2\pi}\int_{0}^{2\pi}\log K_{m}(\tau (s+re^{\sqrt{-1}\theta}))d\theta
\] 
holds. 
Hence $\log K_{m}$ is plurisubharmonic on $\tau (\Delta (s,r))$. 

Next we shall deal with the case that there exists no local section passing 
through $x$. In this case the above argument also holds replacing section by 
holomorphic disk $\tau : \Delta \longrightarrow X$ passing through $x$. 
 
Since $x$ an $\tau$ are arbitrary, the curvature $\Theta_{h_{m}} = \sqrt{-1}\partial\bar{\partial}\log K_{m}$ is semipositive everywhere on $X^{\circ}$.

Since for every local holomorphic section $\sigma$ of $f_{*}{\cal O}_{X}(mK_{X/S})$  the function 
\[
(\sqrt{-1})^{n^{2}}\int_{X_{s}}(\sigma\wedge\bar{\sigma})^{\frac{1}{m}}
\] 
is of algebraic growth along $S - S^{\circ}$.  
More precisely for $s_{0} \in S - S^{\circ}$ as in \cite[p.59 and p. 66]{ka1} there exist positive numbers $C,\alpha ,\beta$ such that 
\begin{equation}\label{eq}
(\sqrt{-1})^{n^{2}}\int_{X_{s}}(\sigma\wedge\bar{\sigma})^{\frac{1}{m}}
\leqq C\cdot \mid s - s_{0}\mid^{-\alpha}\cdot\mid\log (s-s_{0})\mid^{\beta}
\end{equation}
holds. 
Hence $h_{m}$ is also of algebraic growth by its definiton  along the fiber on $S - S^{\circ}$. 
Hence the current $- dd^{c}\log h_{m}$ on $X^{\circ}$ extends to a 
closed positive current on $X$.   

This completes the proof of Theorem 1.1. $\square$ 
\begin{remark}\label{semipos}
As one sees in the proof above, the semipositivity of the curvature $\Theta_{h_{m}}$ is the consequence of the Kawamata's semipositivity (\cite[p.63 Lamma 7]{ka1} and the 
fact that the supremum of plurisubharmonic functions are again plurisubharmonic
if we take the uppersemicontinuous envelope (cf. Theorem \ref{Lelong}).
$\square$ \vspace{5mm}\\
\end{remark}

Let us consider the  direction in which $\Theta_{h_{m}}$ is strictly positive.
Let 
\[
p: \mathbb{P}(F_{m}^{*})\longrightarrow S
\]
be the projective bundle associated with $F_{m}^{*}$ and let
\[
L \longrightarrow \mathbb{P}(F_{m}^{*})
\] 
be the tautological bundle on $L$.
Then 
\[
\parallel\eta \parallel_{\frac{1}{m}} = \mid\!\!\int_{X_{s}}(\eta\wedge\bar{\eta})^{\frac{1}{m}}\!\mid^{\frac{m}{2}} \,\,\,\,\, (\eta\in H^{0}(X_{s},{\cal O}_{X_{s}}(mK_{X_{s}})))
\]
defines a continuous hermitian metric $h_{L}$ on $L$. 
By \cite[p.63, Lemma 7 and p.64, Lemma 8]{ka1} we see that the curvature current 
$\Theta_{h_{L}}$ is semipositive.

We note that on the fiber of $p$, $\Theta_{h_{L}}$ is semipositive and is the curvature  of the tautological line bundle which is ample on the fiber.
In fact let $\eta (t) (t\in \Delta )$ be a holomorphic family of 
sections on a fixed projective variety $X$.
Then by the concavity of logarithm, for every $r\in [0,1)$
\begin{eqnarray*}
\log (\mid\!\!\int_{X}(\eta(0)\wedge\bar{\eta}(0))^{\frac{1}{m}}\!\mid )
& = &\log (\int_{X}(\frac{1}{2\pi}\mid\!\!\int_{0}^{2\pi}(\eta(re^{\sqrt{-1}\theta})\wedge\bar{\eta}(re^{\sqrt{-1}\theta}))^{\frac{1}{m}}\!\mid )d\theta ) \\
& \geqq & 
 \frac{1}{2\pi}\int_{0}^{2\pi}\log (\mid\!\!\int_{X}\eta(re^{\sqrt{-1}\theta})\wedge\bar{\eta}(re^{\sqrt{-1}\theta})^{\frac{1}{m}}\!\mid )d\theta
\end{eqnarray*}
hold  and the similar inequality holds for every point on $S$ (this is a special case of \cite[p. 63, Lemma 7 and p.64, Lemma 8]{ka1}).
 
Hence $\Theta_{h_{L}}\mid_{p^{-1}(s)}$ is strictly positive on the generic point of 
the fiber $p^{-1}(s)$ for every $s\in S$. 
We note that   $\Theta_{h_{L}}$ is strictly positive in the direction in which the coupling of the Kodaira-Spencer class
of the logarithmic deformation in 
the direction and the $m$-th root of the $m$-canonical forms is nonzero 
as in  the proof of \cite[p.63, Lemma 7]{ka1} (this is not explicitly stated 
in the proof of \cite[p.63, Lemma 7]{ka1}, but the discussion in \cite{ka3} clarifies this standard fact).
Hence to specify the strictly posive direction of $\Theta_{h_{m}}$, we need the following infinitesimal Torelli theorem (see also Remark \ref{semipos} above).
\begin{theorem}(\cite{ka3})\label{torelli}
Let $X$ be a normal projective variety of dimension $d$ having only 
canonical singularities. 
We assume that $X$ is a good minimal algegbraic variety, i.e., $m_{0}K_{X}$
 is Cartier and $\hat{\omega}_{X}^{m_{0}}= {\cal O}_{X}(m_{0}K_{X})$
  is generated by global sections for some positive integer $m_{0}$.  
Let $\gamma : X \longrightarrow B$ be the morphism given by the base point free linear system $\mid\!m_{0}K_{X}\!\mid$. 
We take $m_{0}$ large enough so that $B$ is normal and $\gamma$ has connected
fibers. 
Let $C = (\sigma )$ for some $\sigma \in H^{0}(X,\hat{\omega}_{X}^{m})$ 
with $m = am_{0}$ and $a\in \mathbb{N}$.  

Then for $0 < i < m, i\equiv 1$ and if $2i > m+1$, then the kernel of 
the following homomorphism induced by the cup product
\[
\lambda_{X}^{i} :
\mbox{Ext}^{1}_{X}(\Omega_{X}^{1}(\log C),{\cal O}_{X})
\rightarrow \mbox{Hom}(H^{0}(X,\hat{\omega}_{X}^{m+1-i}),
\mbox{Ext}^{1}_{X}(\Omega_{X}^{1}(\log C),\hat{\omega}_{X}^{m+1-i}))
\]
is contained in the image of the natural homomorphism 
\[
H^{1}(B,\gamma_{*}{\cal O}_{X}(T_{X/B}))\rightarrow
\mbox{Ext}^{1}_{X}(\Omega_{X}^{1}(\log C),{\cal O}_{X}).
\]
$\square$
\end{theorem}  
Hence by the proof of Theorem 1.1 and Theorem \ref{torelli}, we have 
the following refined semipositivity theorem. 
 
\begin{theorem}\label{refined semipositivity}
Let $f : X \longrightarrow S$ be a flat projective family of canonical models
of general type 
with only canonical singularities. 
Let $h_{m}$ be the Narashimhan-Simha singular hermitian metric  on $mK_{X/S}$.
Suppose that $mK_{X/S}$ is Cartier and 
$\mid mK_{X/S}\mid$ is relatively very ample.
Then $\Theta_{h_{m}}$ is strictly positive on the generic point of the fiber in the effectively parametrized direction. $\square$ 
\end{theorem}
\begin{remark}
It is not difficult to generalized Theorem \ref{refined semipositivity} 
for all $m\geqq 2$ such that $f_{*}{\cal O}_{X}(mK_{X/S})\neq 0$ by using 
the ring structure of $\oplus_{\ell =0}^{\infty}f_{*}{\cal O}_{X}(\ell K_{X/S})$. 
Similar stirict positivity theorem holds for the flat families of minimal 
projective varieties under the assumption that  minimal model conjecture holds.
$\square$
\end{remark}
\section{Canonical AZD}
In this section, we shall prove that the canonical bundle of  every smooth projective manifold of nonnegtive Kodaira dimension carries a canonical AZD which is functorial 
under any birational morphism.  
Namely we shall prove the following theorem.

\begin{theorem}\label{canAZD}
Let $X$ be a smooth projective variety (not necessary of general type)
of dimension $n$ 
with nonnegative Kodaira dimension. 
Let $K_{m}$ be the m-th Narashimhan-Simha potential for every positive integer $m$.  We set 
\[
C_{m}:= \int_{X}K_{m}^{\frac{1}{m}}.
\]
We define a semipositive $(n,n)$-form $K_{\infty}$ on $X$ by
\[
K_{\infty}(x):= \sup_{m\geqq 1}\frac{1}{C_{m}}K_{m}^{\frac{1}{m}}(x)
\hspace{5mm} (x\in X).
\]
Then $K_{\infty}$ 
exists as a bounded $(n,n)$-form on $X$ ($n= \dim X$) and 
\[
h := \frac{1}{K_{\infty}}
\]
is an AZD of $K_{X}$. $\square$
\end{theorem}
As a corollary, we have the following :
\begin{corollary}\label{invAZD}
Let $X$ be a smooth projective variety with nonnegative Kodaira dimension.
Then there exists an AZD $h$ on $K_{X}$ which depends only on the complex structure of $X$. $\square$
\end{corollary}

\begin{definition}\label{defcanAZD}
Let $X$ be a smooth projective variety with nonnegative Kodaira dimension.
We call the AZD of $K_{X}$ constructed as in Theorem \ref{canAZD} the {\bf canonical AZD} of $K_{X}$. $\square$ 
\end{definition}

\subsection{Construction of the canonical AZD}

The supremum of a family of plurisubharmonic functions 
uniformly bounded from above is known to be again plurisubharmonic, 
if we modify the supremum on a set of measure $0$(i.e., if we take the uppersemicontinuous envelope) by the following theorem of P. Lelong.

\begin{theorem}(\cite[p.26, Theorem 5]{l})\label{Lelong}
Let $\{\varphi_{t}\}_{t\in T}$ be a family of plurisubharmonic functions  
on a domain $\Omega$ 
which is uniformly bounded from above on every compact subset of $\Omega$.
Then $\psi = \sup_{t\in T}\varphi_{t}$ has a minimum 
uppersemicontinuous majorant $\psi^{*}$  which is plurisubharmonic.
We call $\psi^{*}$ the uppersemicontinuous envelope of $\psi$. 
$\square$ \end{theorem}
\begin{remark} In the above theorem the equality 
$\psi = \psi^{*}$ holds outside of a set of measure $0$(cf.\cite[p.29]{l}). 
$\square$ \end{remark}
{\bf Proof of Theorem \ref{canAZD}}.
Since 
\[
\frac{1}{C_{m}}\int_{X}K_{m}^{\frac{1}{m}} = 1
\]
holds for every $m$ such that $\mid mK_{X}\mid\neq\emptyset$. 
Since $h_{m} = 1/K_{m}$ has semipositive curvature in the sense of current, 
by the submeanvalue  property of plurisubharmonic functions, we see that
$\{ C_{m}^{-1}\cdot K_{m}^{\frac{1}{m}}\}_{m=1}^{\infty}$ is a uniformly bounded family of semipositive $(n,n)$-forms on $X$.
Hence the pointwise supremum  
\[  
K_{\infty}:= \sup_{m\geqq 1}\frac{1}{C_{m}}K_{m}^{\frac{1}{m}}
\]
exists as a semipositive $(n,n)$-form on $X$.
By definition
\[
K_{\infty} \geqq \frac{1}{C_{m}}K_{m}^{\frac{1}{m}}
\]
holds for every $m$. 
Hence
\[
h := \mbox{the lower envelope of}\,\,\frac{1}{K_{\infty}}
\]
is a well defined singular hermitian metric on $K_{X}$ and 
\begin{equation}\label{inclusion1}
{\cal I}(h^{m}) \supseteq {\cal I}(h_{m})
\end{equation}
holds for every $m \geqq 1$.
We note that  the definition of $K_{m}$,
\begin{equation}\label{inclusion2}
H^{0}(X,{\cal O}_{X}(mK_{X})\otimes {\cal I}(h_{m}))
= H^{0}(X,{\cal O}_{X}(mK_{X}))
\end{equation}
holds for every $m \geqq 1$. 
Then combining (\ref{inclusion1}) and (\ref{inclusion2}), 
we see that 
\[
H^{0}(X,{\cal O}_{X}(mK_{X})\otimes {\cal I}(h^{m}))
=  H^{0}(X,{\cal O}_{X}(mK_{X}))
\]
holds for every $m\geqq 1$.

On the other hand, by Theorem \ref{Lelong} , we see that 
the curvature current $\Theta_{h}$ of $h$ is seimipositive in the sense of
current. 

Hence $h$ is an AZD of $K_{X}$. $\square$  \vspace{5mm}\\
\subsection{Canonical AZD for families}

The construction of the canonical AZD in Theorem \ref{canAZD} can be applied to
a flat projective family of varieties with only canonical singularities 
similarly as in Section \ref{continuity}.

\begin{theorem}\label{family}
Let $f : X \longrightarrow S$ be  a  flat projective family of varieties 
with connected fibers such that a general fiber of $f$ has  only canonical singularities and  nonnegative Kodaira dimesnion. 
We set  
\[
S^{\circ} := \{ s\in S\mid \mbox{$X_{s}:= f^{-1}(s)$ is reduced with 
only  canonical singularities}\}. 
\]
Then there exists a singular hermitian metric $h$ on $K_{X/S}$ with semipostive curvature  current such that for every  $s\in S^{\circ}$, 
$h\mid_{X_{s}}$ is an AZD of $K_{X_{s}}$. $\square$
\end{theorem}
{\bf Proof of Theorem \ref{family}}.
The construction is essentially same as in the last subsection. 
Let $\{ K_{m,s}\} (s\in S^{\circ})$ be the family of the $m$-th Narashimhan-Simha potentials over $S^{\circ}$. 
We define the relative potential $K_{m}$ on $f^{-1}(S^{\circ})$ by 
\[
K_{m}\!\mid\! X_{s} = K_{m,s} \hspace{5mm}(s\in S^{\circ}).
\]
Suppose that there exists $\sigma_{0} \in H^{0}(X_{s},{\cal O}_{X_{s}}(m_{0}K_{X_{s}}))$ such that 
\[
\mid\int_{X_{s}}(\sigma_{0}\wedge\bar{\sigma}_{0})^{\frac{1}{m_{0}}}\mid = 1
\]
Then for any point $x\in X_{s}$ and positive integer $\ell$, we see that 
\[
K_{m_{0}\ell,s}(x) \geqq \mid\sigma_{0}\mid^{2\ell}(x)
\]
holds by the definition of $K_{m_{0}\ell,s}$.  
Hence by using the invariance of plurigenera (\cite{s1,s2,tu3,ka2}), since we have assumed that every fiber on $S^{\circ}$ has nonnegative Kodaira dimension, there exists a positive  continuous function $c(s)$ on $S^{\circ}$  such that for every positive integer $\ell$ and $s\in S^{\circ}$
\begin{equation}\label{lower}
\int_{X_{s}}K_{m_{0}\ell,s}^{\frac{1}{m_{0}\ell}}\geqq c(s)
\end{equation}
holds.  
 
On the other hand, we obtain the upper estimate of 
\[
\int_{X_{s}}K_{m,s}^{\frac{1}{m}}
\]
as follows. 
Let us fix a $C^{\infty}$-volume form $dV_{s}$ on $X_{s}$, 
such that $dV_{s}^{-1}$ is a $C^{\infty}$ hermitian metric on the $\mathbb{Q}$ line bundle $K_{X_{s}}$. 
Then for any $\sigma \in H^{0}(X_{s},{\cal O}_{X_{s}}(mK_{X_{s}}))$ with 
\[
\mid\int_{X_{s}}(\sigma\wedge\bar{\sigma})^{\frac{1}{m}}\mid = 1,
\] 
by the submeanvalue property of plurisubharmonic functions, 
there exists a positive constant $C^{\prime}(s)$ independent of $\sigma$ and 
$m$ such that  
\[ 
\mid \frac{(\sigma \wedge\bar{\sigma})^{\frac{1}{m}}}{dV_{s}}\mid \leqq C^{\prime}(s)
\]
holds on $X_{s}$.  Hence there exists a positive continuous function  $C(s)$ on $S^{\circ}$ such that 
\begin{equation}\label{upper}
\int_{X_{s}}K_{m,s}^{\frac{1}{m}} \leqq C(s)
\end{equation}
holds for every $m \geqq 1$.

Let $s_{0}\in S^{\circ}$ be an arbitrary point. 
We set  
\[
C_{m} := \int_{X_{s_{0}}}K_{m,s_{0}}^{\frac{1}{m}}. 
\]
Let 
\[
K_{\infty} := \sup_{m\geqq 1} \frac{1}{C_{m}}K_{m}^{\frac{1}{m}}
\]
be the pointwise supremum. 
Then $K_{\infty}^{-1}\mid X_{s_{0}}$ is an AZD on $K_{X_{s_{0}}}$ as in Theorem \ref{canAZD}. 
By the estimates (\ref{lower}) and (\ref{upper}), we see that 
$K_{\infty}$ is a nontrivial locally bounded realative volume form 
on $f^{-1}(S^{\circ})$. 
By the plurisubharmonic variation property of $K_{m}$ (Theorem \ref{main theorem}),  
\[
dd^{c}\log K_{\infty} \geqq 0
\]
holds on $S^{\circ}$. 
Then again by the proof of Theorem \ref{canAZD}, we see that 
$K_{\infty}\mid X_{s}$ is an AZD of $K_{X_{s}}$ for every $s\in S^{\circ}$
by the above construction (especially by the estimates (\ref{lower}) and (\ref{upper})).

Now we consider the extension of the closed positive current 
\[
dd^{c}\log K_{\infty}
\]
across $X - f^{-1}(S^{\circ})$.   
But this is trivial,  since each $K_{m}^{\frac{1}{m}}$ has 
only algebraic singularities along $X - f^{-1}(S^{\circ})$, 
in other words,  for any local holomorphic section $\sigma$ of $f_{*}{\cal O}_{X}(mK_{X/S})$,
\[
(\sqrt{-1})^{n^{2}}\int_{X_{s}}(\sigma \wedge \bar{\sigma})^{\frac{1}{m}}
\hspace{5mm}(n = \dim X - \dim S)
\] 
has algebraic growth glong $S - S^{\circ}$ 
as (\ref{eq}) in Section 3.2  by its definition (this has been already observed in \cite[pp.65-66]{ka1}).
Hence  by  definition, $dd^{c}\log K_{\infty}$ also extends 
across $X - f^{-1}(S^{\circ})$ as a closed positive current (one may also use Lemma \ref{ext} below). 
$\square$  \vspace{5mm} \\ 
By Theorem \ref{main theorem} and the recent result of 
finite generation canonical rings (\cite{b-c-h-m}) 
and Theorem \ref{torelli} we have the following theorem immediately
(see Theorem \ref{refined semipositivity}). 
\begin{theorem}\label{refined semipositivity 2}
Let $f : X \longrightarrow S$ be a flat projective family of varieties of general type with only canonical singularities. 
Let $h$ be the  canonical AZD  on $mK_{X/S}$. 
Then $h$ has semipositive curvature $\Theta_{h}$ in the sense of current 
everywhere on $X$ and $\Theta_{h}$ is strictly positive 
on the generic point of the fiber in the birationally effectively parametrized direction
on $X_{s}$ for every $s\in S_{0}$. $\square$  
\end{theorem}

\subsection{Birational invariance}\label{invariance}

By the definition of the canonical AZD in Theorem \ref{canAZD}, one can see easily that 
it is invariant under birational morphism. 

Let $X$, $X^{\prime}$ be smooth projective variety with nonnegative Kodaira dimension.  
Suppose that $X$ and $X^{\prime}$ are birational. 
Then there exists a smooth projective variety $Y$ which dominates birationally
both $X$ and $X^{\prime}$. 
Let $f : Y \longrightarrow X$ and $f^{\prime} : Y \longrightarrow X^{\prime}$
be the birational morphism. 

For every positive integer $m$, let us denote the $m$-th Narashimhan-Simha potential of $X$ and $X^{\prime}$ by $K_{m}$ and $K_{m}^{\prime}$ respectively.

Then for every positive integer $m$, by the definition of the Narashimhan-Simha
potential, we see that 
\[
f^{*}K_{m} = (f^{\prime})^{*}K_{m}^{\prime}
\]
holds on $Y$. 
Hence the Narashimhan-Simha potential is invariant under birational morphisms
and the canonical AZD for $K_{X}$ is also invariant under birational morphisms.

\subsection{An application of the canonical AZD}
We shall show an easy  application of the canonical AZD. 
\begin{theorem}\label{nef}
Let $f : X \longrightarrow C$ be a semistable projective family 
over a projective curve $C$ such that 
$K_{X/C}\!\mid F$ is ample for every fiber $F$.

Then $K_{X/C}$ is nef. $\square$
\end{theorem}
{\bf Proof of Theorem \ref{nef}}. Let $C_{0}\subset C$ be the Zariski open subset such that
$f$ is smooth over $C_{0}$. 
Let $m$ be a positive integer such that $mK_{X/C}$ is relatively very ample. 
Let $h_{m}$ be the Narashimhan-Simha metric on 
$mK_{X/C}\mid f^{-1}(C_{0})$.  

Let $F_{0} = \sum D_{i}$ be a singular fiber. 
We note that in this case by the $L^{2}$-extension theorem \cite[Theorem]{o}, for every $i$, every element of 
$H^{0}(D_{i},{\cal O}_{D_{i}}(mK_{D_{i}}+ (m-1)\sum_{j\neq i}D_{j}\mid D_{i})$
extends holomorphically to a section 
$H^{0}(f^{-1}(V),{\cal O}_{X}(mK_{X/C}))$, where $V$ is a neighbourhood of 
$f(F_{0})$. 
Hence we see that $h_{m}= K_{m}^{-1}$ extends to a singular hermitian metric on 
$mK_{X/C}$
on the whole $X$ so that $K_{m}$ extends continuously to $X$. 

Then by the definition of $h_{m}$, for every sufficiently large $m$,
$h_{m}\!\!\mid\!\!D_{i}$ is a smooth metric on 
the globally generated sheaf 
${\cal O}_{D_{i}}(mK_{D_{i}}\!+(m-1)\sum_{j\neq i}D_{j}\!\mid D_{i})$. 
By the construction of the canonical AZD $h$, we see that $h$ is a singular hermitian metric on $K_{X/C}$ with the semipositive curvature current $\Theta_{h}$  satisfying 
$\nu (\Theta_{h})\equiv 0$ on $X$  by the relatively ampleness of $K_{X/C}$.
This implies that $K_{X/C}$ is nef by Lemma \ref{nefness}. 
$\square$ 

\section{Positivity of the direct images of pluricanonical systems}  

In this section we prove Theorem \ref{semipositive}.  

\subsection{Theorems of Maitani-Yamaguchi and Berndtsson}

In 2004, Maitani and Yamaguchi proved the following theorem.

\begin{theorem}(\cite{m-y})\label{m-y} Let $\Omega$ be a pseudoconvex domain 
in $\mathbb{C}_{z}\times \mathbb{C}_{w}$ with $C^{1}$ boundary. 
Let $\Omega_{t} := \Omega \cap (\mathbb{C}_{z}\times \{ t\})$ and 
Let $K(z,t)$ be the Bergman kernel function of 
$\Omega_{t}$. 

Then $\log K(z,t)$ is a plurisubharmonic function on $\Omega$. $\square$
\end{theorem}
Recently  generalizing Theorem \ref{m-y}, B. Berndtsson proved the following
higher dimensional and twisted version of Theorem \ref{m-y}.
 
\begin{theorem}(\cite{b1})\label{b1}
Let $D$ be a pseudoconvex domain in $\mathbb{C}^{n}_{z}\times \mathbb{C}^{k}_{t}$.   And let $\phi$ be a plurisubharmonic function on $D$.
For $t\in \Delta$, we set  $D_{t} := \Omega \cap (\mathbb{C}^{n}\times \{ t\})$
and $\phi_{t} := \phi\mid D_{t}$. 
Let $K(z,t) (t\in \mathbb{C}^{k}_{t})$ be the Bergman kernel of the Hilbert space 
\[
A^{2}(D_{t},e^{-\phi_{t}}) := \{ f\in {\cal O}(\Omega_{t})\mid 
\int_{D_{t}}e^{-\phi_{t}}\mid f\mid^{2} < + \infty \} .
\]
Then $\log K(z,t)$ is a plurisubharmonic function on $D$. $\square$
\end{theorem}
As in mensioned in \cite{b2}, his proof also works for a pseudoconvex 
domain in a locally trivial family of manifolds which admits 
a Zariski dense Stein subdomain.  

Also he proved the following theorem. 

\begin{theorem}(\cite[Theorem 1.1]{b2})\label{b2}
Let us consider a domain $D = U \times\Omega$ and let $\phi$ be a 
plurisubharmonic function on $D$.
For simplicity we assume that $\phi$ is smooth up to the boundary and 
strictly plurisubharmonic in $D$. 
Then for each $t\in U$, $\phi_{t} := \phi (\cdot ,t)$ is plurisubharmonic 
on $\Omega$.  Let $A^{2}_{t}$ be the Bergman space of holomorphic functions
on $\Omega$ with norm
\[
\parallel f\parallel^{2} = \parallel f\parallel^{2}_{t}:=  \int_{\Omega} e^{-\phi_{t}}\mid f\mid^{2}. 
\]
The spaces $A^{2}_{t}$ are all equal as vector spaces but have norms 
that vary with $t$. 
Then ``infinite rank'' vector bundle $E$ over $U$ with fiber
$E_{t} = A^{2}_{t}$ is therefore trivial as a bundle but is equipped with 
a notrivial metric. 
Then $(E,\parallel \,\,\,\,\parallel_{t})$ is strictly positive 
in the sense of Nakano. $\square$ 
\end{theorem}
In Theorem \ref{b1} the assumption that $D$ is a pseudoconvex domain 
in the product space is rather strong.  And in Theorem \ref{b2}, Berndtsson
also assumed that $D$ is a product. 

\subsection{Variation of hermitian adjoint bundles}

Recently Berndtsson and I have independently generalized Theorems \ref{b1} and \ref{b2} to the case of projective families. 
For our purpose, we need the following version.  
Theorem \ref{v} has been obtained independently by Berndtsson and Paun \cite{b3,b-p} by
a very different method.

\begin{theorem}\label{v}(\cite{tu7,tu8,b-p})
Let $f : X \longrightarrow S$ be a projective family 
of projective varieties over a complex manifold $S$.
Let $S^{\circ}$  be the maximal nonempty Zariski open subset such that 
$f$ is smooth over $S^{\circ}$.

Let $(L, h)$ be a singular hermitian line bundle on $X$ such that 
$\Theta_{h}$ is semipositive on $X$. \\
Let $K_{s} := K(X_{s},K_{X}+ L\mid_{X_{s}},h\mid_{X_{s}})$ be the Bergman kernel 
of $K_{X_{s}}+ (L\mid X_{s})$ with respect to $h\mid X_{s}$ for $s\in S^{\circ}$. 
Then the singular hermitian metric $h_{B}$ of $K_{X/S} + L\mid f^{-1}(S^{\circ})$ defined by 
\[
h_{B}\mid X_{s}:= K_{s}^{-1} (s\in S^{\circ})
\]
has semipositive curvature on $f^{-1}(S^{\circ})$ and extends on $X$ as a singular hermitian metric on $K_{X/S} +L$ with semipositive curvature current. $\square$ 
\end{theorem}
\begin{theorem}\label{nakano}(\cite{tu7,tu8})
Let $f : X \longrightarrow S$ be a  projective family 
of over a complex manifold $S$
such that $X$ is smooth. 
Let $S^{\circ}$  be a nonempty Zariski open subset such that 
$f$ is smooth over $S^{\circ}$.
Let $(L,h)$ be a hermitian line bundle on $X$ such that 
$\Theta_{h}$ is semipositive on $X$.
We define the hermitian metric $h_{E}$ on  $E:= f_{*}{\cal O}_{X}(K_{X/S} +L)\mid S^{\circ}$
 by 
\[
h_{E}(\sigma ,\tau ) := (\sqrt{-1})^{n^{2}}\int_{X_{s}}h\cdot \sigma\wedge\bar{\tau} \hspace{5mm} (\sigma ,\tau \in H^{0}(X_{s},{\cal O}_{X_{s}}(K_{X_{s}}+L\mid_{X_{s}}))),
\] 
where $n$ denotes the relative dimension of $f : X \longrightarrow S$. 
Let $S_{0}$ be the maximal  Zariski open subset of $S^{\circ}$ such that 
$E\mid S_{0}$ is locally free. 
Then $(E,h_{E})\!\!\mid\!\!S_{0}$ is semipositive in the  sense of Nakano. 
Moreover if $\Theta_{h}$ is strictly positive, then 
$(E,h_{E})\!\!\mid\!S_{0}$ is strictly positive in the sense of Nakano.
$\square$ 
\end{theorem}
Here we reproduce the proofs of Theorems \ref{v} and \ref{nakano} here for the completeness. 
\vspace{3mm}\\
\noindent{\bf  Proof of Theorems \ref{v} and \ref{nakano}}. \vspace{3mm} \\  
Let $f : X \longrightarrow S$ be a projective family. 
Since the statement is local we may assume that $S$ is the unit open ball 
$B$ with center $O$ in {$\mathbb{C}^{m}$. 
We may also and do assume that the family $f : X \longrightarrow B$ is 
a restriction of a projective family 
\[
\hat{f} : \hat{X} \longrightarrow B(O,2)
\]
over the open ball $B(O,2)$ of radius $2$ with center $O$.
Let 
\[
F : \hat{X} \times B(O,2) \longrightarrow B(O,2)\times B(O,2)
\]
the fiber space defined by 
\[
F(x,t) = (f(x),t). 
\]
Let $\varepsilon$ be a positive number less than $1$.
We set   
\[
T(\varepsilon ) = \{ (s,t) \in B(O,2)\times B(O,2) \mid 
t\in B, s\in B(t,\varepsilon )\}
\]
and  
\[
X(\varepsilon ) := F^{-1}(T(\varepsilon )).
\]
Let 
\[
f_{\varepsilon} : X(\varepsilon ) \longrightarrow B
\]
be the family  defined by 
\[
f_{\varepsilon}(x,t) = t.
\]
Since for $(x,t)\in X(\varepsilon )$, $x\in f^{-1}(B(t,\varepsilon ))$ holds,
we see that  
\[
X(\varepsilon )_{t} := f^{-1}(B(t,\varepsilon ))
\]
holds.  Hence we may consider $X(\varepsilon )_{t}$ as a family 
\[
\pi_{\varepsilon ,t}: X(\varepsilon ,t)\longrightarrow B(t,\varepsilon ).
\] 
We note that $T(\varepsilon )$ is a domain of holomorphy 
in $\mathbb{C}^{2m}$.  Hence $X(\varepsilon )$ is a pseudoconvex domain in $X\times B(O,2)$.  Since $X \times B(O,2)$ is a product manifold,  the proof of Theorem \ref{b1} works without any essential change in this case (cf. \cite{b1}).
Hence if we define $K_{\varepsilon}$ by 
\[
K_{\varepsilon}\mid X(\varepsilon )_{t}:= K(X(\varepsilon )_{t},K_{X}+L\mid X(\varepsilon )_{t},h_{L}\mid X(\varepsilon )_{t}), 
\]
then 
\begin{equation}\label{psh}
\sqrt{-1}\partial\bar{\partial}\log K_{\varepsilon}
 \geqq 0
\end{equation}
holds on $X(\varepsilon )$. 
We note that 
\begin{equation}\label{limit}
\lim_{\varepsilon\downarrow 0}\mbox{vol}(B(t,\varepsilon ))\cdot K(X(\varepsilon )_{t},K_{X}+L\mid\!X(\varepsilon )_{t},h_{L}\mid\!X(\varepsilon )_{t})
= K(X_{t},K_{X_{t}} + L\mid X_{t}, h\mid X_{t})\wedge f^{*}d\mu (t)
\end{equation}
holds for almost all $t\in B$, where $d\mu (t)$ denotes the standard volume form on $\mathbb{C}^{m}$. 
In fact, if we consider the family 
\[
\pi_{\varepsilon ,t} : X(\varepsilon ) \longrightarrow B(t,\varepsilon )
\] 
as a family over the unit open ball  $B$ in $\mathbb{C}^{m}$ with center $O$  by 
\[
t^{\prime}\mapsto \varepsilon^{-1}(t^{\prime} - t),
\]
the limit as $\varepsilon \downarrow 0$ is nothing but the 
trivial family $X_{t}\times B$.   
And we note that  for a $L$-valued canonical form $\sigma$ on $f^{-1}(B(t,\varepsilon ))$,
\[
\int_{B(t,\varepsilon )}h_{E}(\sigma,\sigma)
\]
is nothing but the $L^{2}$-norm of the $L$-valued canonical form $\sigma$ with respect $h_{L}$ over $f^{-1}(B(t,\varepsilon ))$ by Fubini's theorem, 
where we abbrebiate the standard Lebesgue measure on $\mathbb{C}^{m}$.
In this way the $L^{2}$-norm on  $f^{-1}(B(t,\varepsilon ))$ can be viewd 
the molification of the $L^{2}$-norm on the fiber $f^{-1}(t)$.  
Hence as $\varepsilon\downarrow 0$, the molified norm converges to the 
norm on the fiber almost everywhere on $B$.  
 Then the desired equality  follows from 
the $L^{2}$-extension theorm (\cite{o-t,o}) and the extremal property of the 
Bergman kernels.  In fact the $L^{2}$-extension theorem implies 
\[
\lim_{\varepsilon\downarrow 0}\mbox{vol}(B(t,\varepsilon ))\cdot K(X(\varepsilon )_{t},K_{X}+L\mid\!X(\varepsilon )_{t},h_{L}\mid\!X(\varepsilon )_{t})
\geqq  K(X_{t},K_{X_{t}} + L\mid X_{t}, h\mid X_{t})\wedge f^{*}d\mu (t)
\]
for almost all $t\in B$ and the converse inequality follows from the convergence to the trivial family as above.

Combining (\ref{psh}) and (\ref{limit}), we see that 
$h_{B}$ has semipositive curvature on $f^{-1}(S^{\circ})$. 

To get an extension of $h_{B}$ across $X - f^{-1}(S^{\circ})$, let us 	quote the following lemma. 

\begin{lemma}(\cite[Corollary 7.3]{b-t})\label{ext}
Let $\{ u_{j}\}$ be a sequence of plurisubharmonic functions locally bounded above 
on the bounded open set $\Omega$ in $\mathbb{C}^{m}$. Suppose further 
\[
\limsup_{j\rightarrow\infty}u_{j}
\]
is not identically $-\infty$ on any component of $\Omega$. 
Then there exists a plurisubharmonic function $u$ on $\Omega$
such that the set of points
\[
\{x \in \Omega \mid u(x) \neq (\limsup_{j\rightarrow\infty}u_{j})(x)\}
\]
is pluripolar.
$\square$
\end{lemma}

\noindent 

The extension of $h_{B}$ in Theorem \ref{v} follows from Lemma \ref{ext}. In fact since \\
$K(X(\varepsilon )_{t},K_{X}+L\mid\!X(\varepsilon )_{t},h_{L}\mid\!X(\varepsilon )_{t})$ varies plurisubharmonic way in $t\in B$, 
\[
\{ \mbox{vol}(B(t,\varepsilon ))\cdot K(X(\varepsilon )_{t},K_{X}+L\mid\!X(\varepsilon )_{t},h_{L}\mid\!X(\varepsilon )_{t})\mid_{X_{t}}\}_{t\in B}
\]
is uniformly bounded on $B$ even when $t \in S - S^{\circ}$. 
Then  letting $\varepsilon$ tend to $0$ and applying Lemma \ref{ext} 
we see that $-dd^{c}\log h_{B}$ extends across $X - f^{-1}(S^{\circ})$
as a closed positive current. 

This completes the proof of  Theorem \ref{v}. \vspace{3mm}

The proof of Theorem \ref{nakano}, is quite similar. 
First we note that 
\[
f_{\varepsilon} : X(\varepsilon ) \longrightarrow B
\]
is everywhere smooth.  

For the moment,  we shall assume that $E$ is locally free on $B(O,2)$.
Then there exists a  global generator $\{ \sigma_{1},\cdots ,\sigma_{r}\}$
of $E$ on $B(O,2)$, where $r = \mbox{rank}\, E$. 
Then we see that every $t\in B$, the fiber of the vector bundle 
$(f_{\varepsilon})_{*}{\cal O}_{X(\varepsilon )}(K_{X(\varepsilon )/B}+L)\otimes{\cal I}(h_{L})$ at $t$ is  canonically isomorphic to 
$\mathbb{C}^{r}\times {\cal O}(B(t,\varepsilon ))$ in terms of the frame
$\{ \sigma_{1},\cdots ,\sigma_{r}\}$.  And moreover the space
${\cal O}(B(t,\varepsilon ))$ is canonically isomorphic to ${\cal O}(B(O,\varepsilon ))$ by the parallel translation.  
In this case 
\[
E_{\varepsilon}: = (f_{\varepsilon,(2)})_{*}{\cal O}_{X(\varepsilon)}(K_{X(\varepsilon )/B}\otimes p_{1}^{*}L\otimes {\cal I}(p_{1}^{*}h))
\]
is a vector bundle of infinite rank on $B$, where 
\[
p_{1}: X(\varepsilon ) \longrightarrow X
\]
denotes  the first projection
\[
p_{1}(x,t) = x,\hspace{5mm}(x,t)\in  X(\varepsilon )
\]
and $(f_{\varepsilon,(2)})_{*}{\cal O}_{X(\varepsilon)}(K_{X(\varepsilon )/B}\otimes p_{1}^{*}L\otimes {\cal I}(p_{1}^{*}h))$ denotes the direct image 
of $L^{2}$-holomorphic sections. 
 
By the same proof as Theorem \ref{b2}, we see that the curvature current of 
$h_{E,\varepsilon}$ is well defined everywhere on $B$ and is semipositive
in the sense of Nakano.  
Letting $\varepsilon$ tend to $0$, the curvature $\Theta_{h_{E_{\varepsilon}}}$ converges to the curvature $\Theta_{h_{E}}$ operating 
on $E_{t}\otimes {\cal O}_{B,t}$ in the obvious manner
for every $t$ such that $f$ is smooth over $t$. 
Hence $\Theta_{h_{E}}$ is semipositive in the sense of Nakano. 

This completes the proof of Theorems \ref{v} and \ref{nakano}. $\square$ 
\subsection{Proof of Theorem \ref{semipositive}}

Let $f : X \longrightarrow S$ be a flat projective family with connected fibers  such that a general fiber of $f$ has only canonical singularities. 
We assume $X$ is normal and $S$ is smooth.  Let $n$ denote the relative dimension of $f: X \longrightarrow S$.
We set  
\[
S^{\circ} := \{ s\in S\mid \mbox{$X_{s} := f^{-1}(s)$ is reduced and has only canonical singularities}\}.
\] 
Then $K_{X/S}$ has a canonial AZD $h$ over $S$ such that 
$\Theta_{h}$ is semipositive on $X$ (Theorem \ref{family}). 
And over $S^{\circ}$, $F_{m}$ is locally free by \cite{ka2,s2, tu3}.

Then 
the hermitian inner product $h_{F_{m}}$ defined by  
\[
h_{F_{m}}(\sigma ,\sigma^{\prime})_{s}
:= (\sqrt{-1})^{n^{2}}\int_{X_{s}}\sigma\wedge\bar{\sigma}^{\prime}\cdot h^{m-1}
\]
is continuous  on $S^{\circ}$, because $h^{m}(\sigma ,\sigma^{\prime})$ is bounded on 
$X_{s}$ by the construction of $h$ (cf. Proof of Theorem \ref{family}).  
And the curvature current of $h_{F_{m}}$ is semipositive in the sense of  Nakano on  $S^{\circ}$ by Theorem \ref{nakano} (here $X$ need not be smooth but Theorem \ref{nakano} is applicable by resolution and the continuity of metrics (\cite{ka2,s2, tu3})).
This completes the proof of Theorem \ref{semipositive}. $\square$ 
\begin{remark}
In Theorem \ref{semipositive}, we have only considered the relative 
canonical AZD $h$ constructed in Theorem \ref{family}. 
But of course the smilar statement holds for  the $L^{2}$-inner product:
\[
h^{\prime}_{F_{m}}(\sigma ,\sigma^{\prime})_{s}: = 
(\sqrt{-1})^{n^{2}}\int_{X_{s}}\sigma\wedge \bar{\sigma}^{\prime}\cdot h_{m}^{\frac{m-1}{m}}, 
\]
where $h_{m}$ is as in Theorem \ref{main theorem}.
This may be useful, since $h_{m}$ has only algebraic singularities. 
One may also consider many other variants. 
See also Remark \ref{fg}. $\square$
\end{remark}

Author's address\\
Hajime Tsuji\\
Department of Mathematics\\
Sophia University\\
7-1 Kioicho, Chiyoda-ku 102-8554\\
Japan \\
e-mail address: tsuji@mm.sophia.ac.jp

\end{document}